\documentclass[reqno]{amsart}
\usepackage{a4wide}
\usepackage{amsmath, amsthm, amssymb}
\usepackage[margin=2cm]{geometry}
\usepackage{graphicx}
\usepackage{hyperref,accents,mathtools}
\usepackage{multicol,caption}
\usepackage{graphicx}
\usepackage{esint}
\usepackage{enumitem}
\newtheorem{theorem}{Theorem}[section]

\newtheorem{lemma}[theorem]{Lemma}
\newtheorem{proposition}[theorem]{Proposition}

\theoremstyle{remark}
\newtheorem{remark}{Remark}[section]

\theoremstyle{definition}
\newtheorem{definition}{Definition}[section]

 \newcommand{\R}{{\mathbb{R}}}
 \newcommand{\N}{{\mathbb{N}}}
 \newcommand{\C}[1]{\mathbf{C^{#1}}}
 \renewcommand{\L}[1]{\mathbf{L^{#1}}}
 \newcommand{\lip}{\mathrm{Lip}}
 \newcommand{\Cc}[1]{\mathbf{C_c^{#1}}}
 \renewcommand{\d}{{\rm{d}}}
 \newcommand{\modulo}[1]{{\left|#1\right|}}
 \newcommand{\BV}{\mathbf{BV}}
 \newcommand{\tv}{\mathrm{TV}}
 \newcommand{\sign}{\mathrm{sign}}
 \newcommand{\supp}{\mathrm{supp}}
 \newcommand{\Lloc}[1]{\mathbf{L^{#1}_{loc}}}
 \newcommand{\BVloc}{\mathbf{BV_{loc}}}

\begin{document}

\title{Deterministic particle approximation of scalar conservation laws.
}


\author{M.Di Francesco, S.Fagioli, and M.D.Rosini}

 \email{marco.difrancesco@univaq.it}
 \email{simone.fagioli@dm.univaq.it}
 \email{mrosini@umcs.lublin.pl}


%\address{Marco Di Francesco.
%              DISIM, Universit\`a degli Studi dell'Aquila,
%              via Vetoio~1 (Coppito), 67100 L’Aquila (AQ), Italy \\
%%              Tel.: +123-45-678910\\
%%              Fax: +123-45-678910\\
%              \email{marco.difrancesco@univaq.it}
%           Simone Fagioli.
%              DISIM, Universit\`a degli Studi dell'Aquila,
%              via Vetoio~1 (Coppito), 67100 L’Aquila (AQ), Italy \\
%%              Tel.: +123-45-678910\\
%%              Fax: +123-45-678910\\
%              \email{simone.fagioli@univaq.it}
%           Massimiliano D. Rosini.
%              Instytut Matematyki, Uniwersytet Marii Curie-Sk\l odowskiej,
%              pl.\ Marii Curie-Sk\l odowskiej~1, 20-031 Lublin, Poland \\
%%              Tel.: +123-45-678910\\
%%              Fax: +123-45-678910\\
%              \email{mrosini@umcs.lublin.pl}
%}
\begin{abstract}
In this paper we prove that the unique entropy solution to a scalar nonlinear conservation law with strictly monotone velocity and nonnegative initial condition can be rigorously obtained as the large particle limit of a microscopic follow-the-leader type model, which is interpreted as the discrete Lagrangian approximation of the nonlinear scalar conservation law.
The result is complemented with some numerical simulations.
\keywords{Scalar conservation laws \and Follow-the-leader system \and Particle approximation}
\end{abstract}

\maketitle

\section{Introduction}

The approximation of scalar nonlinear conservation laws
\begin{equation}\label{eq:CLintro}
  \rho_t + (\rho \, v(\rho))_x = 0
\end{equation}
via microscopic modeling is a longstanding challenge.
A probabilistic approach to this problem has been proposed in a vast literature in the past decades, see e.g.~\cite{ferrari,ferrari_TASEP,landim} and the references therein.
The kinetic approximation of nonlinear conservation laws has been carried out in \cite{lions_perthame_tadmor}.

In \cite{DF_rosini}, the microscopic Lagrangian formulation of \eqref{eq:CLintro} via the follow-the-leader particle system
\begin{equation}\label{eq:FTLintro}
  \dot{x}_i= v\left(\frac{\ell}{x_{i+1}-x_i}\right)
\end{equation}
has been \emph{rigorously} derived for the first time under the assumption that $v$ is monotone decreasing (plus some additional assumptions, see \eqref{V1} and \eqref{V2}).
The derivation is restricted to nonnegative, bounded, and compactly supported solutions $\rho$.
Roughly speaking, the main result in \cite{DF_rosini} states what follows.
Let $\bar\rho \in \L\infty(\R;\R_+)$ be compactly supported. Assume for simplicity that $\bar\rho$ has unit mass. For a given integer $n \in \N$ sufficiently large, let the minimal interval $[\bar{x}_{\min},\bar{x}_{\max}]$ containing $\supp[\bar\rho]$ be split into $n$ intervals containing the mass $\ell_n \doteq 1/n$.
%being $M \doteq \|\bar\rho\|_{\L1(\R)}$ the total mass of $\bar\rho$.
Let the edges of those intervals $\bar{x}_0 \doteq \bar{x}_{\min} < \bar{x}_1 < \ldots < \bar{x}_{n-1} < \bar{x}_n \doteq \bar{x}_{\max}$ be the initial positions of a set of particles with equal mass $\ell_n$.
Let the particles $x_0(t),\ldots,x_{n-1}(t)$ evolve via \eqref{eq:FTLintro} with $\ell=\ell_n$, and let $x_n(t)=\bar{x}_n + v(0) \, t$.
Then, the discretised density
\[
\rho^n(t,x) \doteq \sum_{i=0}^{n-1}\frac{\ell_n}{x_{i+1}(t)-x_i(t)} ~ \mathbf{1}_{[x_i(t),x_{i+1}(t))}
\]
converges up to a subsequence a.e.~in $\Lloc1(\R_+\times \R)$ to the unique entropy solution $\rho$ to \eqref{eq:CLintro} with initial condition $\bar\rho$, see Definition~\ref{def:entropy_sol} below.
Moreover, the empirical measure
\[\tilde{\rho}^n(t) \doteq \sum_{i=0}^{n-1} \ell_n \, \delta_{x_i(t)}\]
converges to $\rho$ in $\Lloc1(\R_+\,;\, d_{1})$, where $d_{1}$ is the $1$-Wasserstein distance on $\R$.

%This note aims at
%\begin{itemize}[itemindent=*,leftmargin=0pt]\setlength{\itemsep}{0cm}%
%  \item shortening the proof of the result in \cite{DF_rosini}, in particular by avoiding the Eulerian-to-Lagrange coordinates change of variables,
%  \item removing the assumption of initial compact support,
%  \item complementing the results of \cite{DF_rosini} with some numerical simulations.
%\end{itemize}
This note aims at shortening the proof of the result in \cite{DF_rosini} (in particular by avoiding the Eulerian-to-Lagrangian coordinates change of variables), removing the assumption of initial compact support and complementing the results of \cite{DF_rosini} with some numerical simulations.

\section{Preliminaries and result}%\label{sec:prelim}

Let us consider the Cauchy problem for a one-dimensional scalar conservation law
\begin{equation}\label{eq:cauchy}
\begin{cases}
  \rho_t + f(\rho)_x = 0,&  (t,x)\in(0,+\infty)\times\R,\\
  \rho(0,x)=\bar\rho(x) & x\in\R,
\end{cases}
\end{equation}
where  $f(\rho) \doteq \rho \, v(\rho)$.
The initial datum $\bar\rho$ and the velocity map $v \colon \R_+ \rightarrow \R$ satisfy the basic assumptions
\begin{align}\tag{I1}\label{I1}
&\bar\rho \in \L\infty(\R)\cap \L1(\R),&
&\bar\rho\geq 0,
\\
\tag{V1}\label{V1}
&v\in \C1(\R_+),&
&v'(\rho)<0 \text{ for } \rho>0.
\end{align}
In some cases, we require the additional (optional) assumptions
\begin{align}\tag{I2}\label{I2}
&\bar\rho\in \BV(\R),
\\
&\tag{V2}\label{V2}
\R_+ \ni \rho \mapsto [\rho \, v'(\rho)] \in \R_- \text{ is non-increasing}.
\end{align}
For simplicity, we shall normalise the total mass and assume $\|\bar\rho\|_{\L1(\R)} = 1$.
%\[\int_\R\bar\rho(x) \d x = 1.\]
We introduce the notation $v_{\max} \doteq v(0)$ and we shall assume for simplicity that $v_{\max}>0$.

%\begin{definition}\label{def:entropy_sol}
%Let $\bar\rho$ satisfy \eqref{I1}. A function $\rho\in \L\infty(\R_+\times\R)$ is called \emph{entropy solution} to the Cauchy problem \eqref{eq:cauchy} if
%\[
%  \iint_{\R_+\times\R} \, \Bigl[|\rho(t,x)-k| \, \varphi_t(t,x) + \sign(\rho(t,x) - k) \bigl[ f(\rho(t,x)) - f(k)\bigr] \varphi_x(t,x)\Bigr] \d x\d t\geq 0
%\]
%for all $\varphi\in \Cc\infty((0,+\infty) \times \R)$ with $\varphi\geq 0$ and for all $k\geq 0$, and additionally
%\[%\begin{equation}\label{eq:initial_condition}
% \lim_{t \searrow 0} \int_\R\rho(t,x) \, \psi(x) \d x = \int_\R\bar\rho(x) \, \psi(x) \d x
%\]%\end{equation}
%for all $\psi\in \Cc\infty(\R)$.
%\end{definition}
\begin{definition}\label{def:entropy_sol}
Let $\bar\rho$ satisfy \eqref{I1}. A weak solution $\rho\in \L\infty(\R_+\times\R)$ to the Cauchy problem \eqref{eq:cauchy} is called \emph{entropy solution} to the Cauchy problem \eqref{eq:cauchy} if
\[
  \iint_{\R_+\times\R} \, \Bigl[|\rho(t,x)-k| \, \varphi_t(t,x) + \sign(\rho(t,x) - k) \bigl[ f(\rho(t,x)) - f(k)\bigr] \varphi_x(t,x)\Bigr] \d x\d t\geq 0
\]
for all $\varphi\in \Cc\infty((0,+\infty) \times \R)$ with $\varphi\geq 0$ and for all $k\geq 0$.
\end{definition}
We point out that the above definition is slightly weaker than the %original
definition in \cite{kruzkov}.
The next theorem collects the uniqueness result in \cite{kruzkov} and its variant in \cite{chen_rascle}.
\begin{theorem}[\cite{chen_rascle,kruzkov}]%\label{thm:uniqueness}
Assume that \eqref{I1} and \eqref{V1} are satisfied.
Then there exists a unique entropy solution according to Definition~\ref{def:entropy_sol}.
\end{theorem}

We now introduce the approximation scheme.
For future use, we introduce the notation
\[R \doteq \|\bar\rho\|_{\L\infty(\R)}.\]
For a given $n\in\mathbb{N}$ sufficiently large, we set $\ell_n \doteq 1/n$.
Let $\bar{x}^n_1$ be defined by
\[\bar{x}^n_1 \doteq \sup\left\{x\in\R\colon\int_{-\infty}^x\bar\rho(x) \d x<\ell_n\right\},\]
and the points $\bar{x}^n_i$ with $i \in \{2,\ldots,n-1\}$ be defined recursively by
\[\bar{x}^n_i = \sup\left\{x\in \R\colon\int_{\bar{x}^n_{i-1}}^x\bar\rho(x) \d x<\ell_n\right\}.\]
It follows that $\bar{x}^n_1 < \bar{x}^n_2 < \ldots < \bar{x}^n_{n-1}$.
Moreover
\begin{align}\label{eq:mass}
&\int_{-\infty}^{\bar{x}^n_1}\bar\rho(x) \d x = \int_{\bar{x}^n_{i-1}}^{\bar{x}^n_{i}}\bar\rho(x) \d x = \int_{\bar{x}^n_{n-1}}^{+\infty}\bar\rho(x) \d x = \ell_n,
&i \in \{2,\ldots,n-1\}.
\end{align}
We let the $(n-1)$ particles defined above evolve according to the follow-the-leader system of ODEs
\begin{align}\label{eq:FTL}
 &\begin{cases}
    \dot{x}_i^n(t)=v(R^n_i(t)), & i \in \{1,\ldots,n-2\}, \\
    \displaystyle{\dot{x}_{n-1}^n(t)=v_{\max}},  \\
    x^n_i(0) = \bar{x}^n_i, & i \in \{1,\ldots,n-1\},
  \end{cases}
  &R^n_i(t) \doteq \frac{\ell_n}{x^n_{i+1}(t)-x^n_i(t)}.
\end{align}
The discrete maximum principle in \cite[Lemma~1]{DF_rosini} ensures the solution $(x_i^n)_{i=1}^{n-1}$ to \eqref{eq:FTL} is well defined, since the particles $(x^n_i)_{i=1}^{n-1}$ strictly preserve their initial order. More precisely, we have the following lemma.
\begin{lemma}[Discrete maximum principle \cite{DF_rosini}]\label{lem:maximum}
  Assume (I1) and (V1) are satisfied. Then, for all $t \in \R_+$, the solution to \eqref{eq:FTL} satisfies
  \begin{align*}
  &x^n_{i+1}(t)-x^n_i(t)\geq \frac{\ell_n}{R},
  &i \in \{1,\ldots, n-2\}.
  \end{align*}
\end{lemma}

We have split the initial condition into $n$ regions with equal mass $\ell_n$. We have then defined the motion of $(n-1)$ particles. This permits to reconstruct a time-depending (piecewise constant) density within the interval $[x^n_1(t),x^n_{n-1}(t)]$, which will consist of $(n-2)$ constant values on as many intervals. Under the natural assumption that a mass $\ell_n$ will be maintained on each interval, we still need to assign mass to two points outside the interval $[x^n_1(t),x^n_{n-1}(t)]$ in order to obtain a time-depending density with unit mass. To perform this task, we set two \emph{artificial particles} $x_0^n(t)$ and $x_n^n(t)$ as follows
\begin{align}\label{eq:phantom}
  &x^n_0(t) \doteq 2x^n_{1}(t)-x^n_2(t),
  &x^n_n(t) \doteq 2x^n_{n-1}(t)-x^n_{n-2}(t),
\end{align}
and let $R^n_0(t) \doteq R^n_1(t)$ and $R^n_{n-1}(t) \doteq R^n_{n-2}(t)$ for all $t\geq 0$.
We then set
\begin{equation}\label{eq:disdens}
\rho^n(t,x) \doteq \sum_{i=0}^{n-1}R^n_i(t) ~ \mathbf{1}_{[x^n_{i}(t),x_{i+1}^n(t))}(x) = \sum_{i=0}^{n-1} \frac{\ell_n}{x^n_{i+1}(t)-x^n_i(t)} ~ \mathbf{1}_{[x^n_{i}(t),x_{i+1}^n(t))}(x).
\end{equation}
We notice that $\int_\R \rho^n(t,x) \d x = n \, \ell_n =1$ and that $\rho^n(t,\cdot)$ is compactly supported for all $n$ and for all $t$.
%Moreover, due to \eqref{eq:phantom}, we can easily see that $R^n_0(t)=R^n_1(t)$ and $R^n_{n-1}(t)=R^n_{n-2}(t)$ for all $t\geq 0$.
For future use we compute
\begin{equation}\label{eq:ODE_density}
  \begin{cases}
  \dot{R}^n_i(t)=-\dfrac{R^n_i(t)^2}{\ell_n} \, \bigl[v(R^n_{i+1}(t))-v(R^n_i(t))\bigr], &i \in \{1,\ldots,n-3\}, \\[10pt]
  \dot{R}^n_{n-2}(t)=-\dfrac{R^n_{n-2}(t)^2}{\ell_n} \, \bigl[v_{\max}-v(R^n_{n-2}(t))\bigr].
  \end{cases}
\end{equation}

\begin{remark}\label{rem01}
In case $\supp[\bar\rho]$ is bounded either from above or from below, it is possible to improve the above construction.
In the former case, the particle $x^n_n$ can be set on
%$\max \{\overline{\supp[\bar\rho]}\}$
$\max \{\supp[\bar\rho]\}$ initially and let evolve with maximum speed $v_{\max}$, and the preceding particle $x^n_{n-1}$ let evolve according to $\dot{x}^n_{n-1}(t)=v(\ell_n/(x^n_n(t)-x^n_{n-1}(t)))$.
In the latter case, the particle $x^n_0$ can be set on
%$\min \{\overline{\supp[\bar\rho]}\}$
$\min \{\supp[\bar\rho]\}$ initially and let evolve according to $\dot{x}^n_{0}(t)=v(\ell_n/(x^n_1(t)-x^n_{0}(t)))$.
In \cite{DF_rosini} both these conditions are required for the initial datum and such construction is applied.
\end{remark}

Our result, which extends the one in \cite{DF_rosini}, reads as follows.
\begin{theorem}\label{thm:main}
%Assume $\bar\rho$ satisfies \eqref{I1} and $v$ satisfies \eqref{V1}.
%Assume also either $\bar\rho$ satisfies also \eqref{I2}, or $v$ satisfies also \eqref{V2}.
Assume that \eqref{I1} and \eqref{V1} are satisfied.
Moreover, assume that at least one of the two conditions \eqref{I2} and \eqref{V2} is also satisfied.
Then, $\rho^n$ converges (up to a subsequence) almost everywhere and in $\Lloc1$ on $\R_+\times \R$ to the unique entropy solution $\rho$ to the Cauchy problem \eqref{eq:cauchy} according to Definition~\ref{def:entropy_sol}.
\end{theorem}
The result in \cite{DF_rosini} also states the convergence of the empirical measure $\tilde{\rho}^n(t) \doteq \sum_{i=1}^{n} \ell_n \, \delta_{x_i^n(t)}$ towards the entropy solution $\rho$.
For the sake of brevity, we shall skip that part in this note.

\section{Proof of the main result}%\label{sec:proof}

In this section we prove Theorem~\ref{thm:main}.
Clearly, the result in Lemma~\ref{lem:maximum} ensures that $\|\rho^n(t,\cdot)\|_{\L\infty(\R)}\leq R$ for all $t\geq0$.
For notational simplicity, whenever it is clear from the context, we shall omit the $n$-dependence in the approximating scheme.
Moreover, as our results is a slight extension of the one in \cite{DF_rosini}, we shall often shorten proofs and refer to the corresponding results in \cite{DF_rosini}, still trying to keep this note as much self-contained as possible.

As usual in the context of scalar conservation laws, a uniform control of the $\BV$ norm is necessary in order to gain enough compactness of the approximating scheme.
In our case, the compactness can be obtained in two distinct ways.
The first one is a uniform $\BV$ contraction property for $\rho^n$, and it obviously requires $\BV$ initial data.

\begin{proposition}\label{pro:bv_contraction}
%  Assume $\bar\rho$ satisfies \eqref{I1}, \eqref{I2} and $v$ satisfies \eqref{V1}.
  Assume that \eqref{I1}, \eqref{I2} and \eqref{V1} are satisfied.
  Then, for all $n\in\mathbb{N}$ one has
  \[\tv[\rho^n(t,\cdot)]\leq \tv[\rho^n(0,\cdot)]\leq \tv [\bar\rho].\]
\end{proposition}

\begin{proof}
The estimate $\tv[\rho^n(0,\cdot)]\leq \tv [\bar\rho]$ is a simple exercise. We now compute
\begin{align*}
  &~\frac{d}{dt}\tv[\rho^n(t,\cdot)]
  = \frac{d}{dt} \left[ R_1(t)+R_{n-2}(t)+\sum_{i=1}^{n-3}|R_i(t)-R_{i+1}(t)|\right]\\
  =&~ \dot{R}_1(t)+\dot{R}_{n-2}(t)+ \sum_{i=1}^{n-3} \sign\bigl(R_i(t)-R_{i+1}(t)\bigr) \bigl[\dot{R}_i(t)-\dot{R}_{i+1}(t)\bigr]\\
  =&~ \bigl[ 1+\sign\bigl(R_1(t)-R_2(t)\bigr) \bigr] \dot{R}_1(t)
  + \bigl[ 1-\sign\bigl(R_{n-3}(t)-R_{n-2}(t)\bigr) \bigr] \dot{R}_{n-2}(t)\\
  &~ + \sum_{i=2}^{n-3} \, \bigl[ \sign\bigl(R_i(t)-R_{i+1}(t)\bigr) - \sign\bigl(R_{i-1}(t)-R_i(t)\bigr) \bigr] \dot{R}_i(t).
\end{align*}
By plugging \eqref{eq:ODE_density} into the above computation and employing the assumption \eqref{V1} one can easily prove that the above quantity is not positive.
\end{proof}

The second way to achieve compactness is via the following \emph{discrete Oleinik-type inequality}. Here we do not require the extra assumption \eqref{I2} on the initial condition, but we need the assumption \eqref{V2} on the velocity map.

\begin{proposition}\label{pro:oleinik}
%  Assume $\bar\rho$ satisfies \eqref{I1} and $v$ satisfies \eqref{V1}, \eqref{V2}.
  Assume that \eqref{I1}, \eqref{V1} and \eqref{V2} are satisfied.
  Then, for all $t\geq0$ one has% the following inequality
  \begin{align}\label{eq:oleinik}
    &\frac{\dot{x}^n_{i+1}(t)-\dot{x}^n_i(t)}{x^n_{i+1}(t)-x^n_i(t)}\leq \frac{1}{t},
    &i\in\{0,\ldots,n-1\}.
  \end{align}
%  \begin{align}\label{eq:oleinik}
%    &\frac{v(R_{i+1}(t))-v(R_i(t))}{x_{i+1}(t)-x_i(t)}\leq \frac{1}{t},
%    &i=0,\ldots,n-2.
%  \end{align}
\end{proposition}

\begin{proof}
Due to \eqref{eq:phantom}, it suffices to prove \eqref{eq:oleinik} for $i\in\{1,\ldots,n-2\}$.
We start by observing that this is equivalent to prove
\begin{align*}%\label{eq:oleinik2}
&z_i(t) \doteq t \, R_i(t) \Bigl[\dot{x}_{i+1}(t)-\dot{x}_i(t)\Bigr] \leq \ell_n,
&i\in\{1,\ldots,n-2\}.
\end{align*}
We shall prove the above estimate inductively on $i$ by using the equations \eqref{eq:ODE_density}.
We drop the time dependency for simplicity.

We start by proving $z_{n-2}= t \, R_{n-2} [ v_{\max}-v(R_{n-2})] \leq \ell_n$.
We have, due to \eqref{eq:ODE_density} and \eqref{V1}, that
\begin{align*}
  &~ \dot{z}_{n-2}
  =   R_{n-2} \, \bigl[v_{\max}-v(R_{n-2})\bigr] + t \, \dot{R}_{n-2} \, \bigl[ v_{\max}-v(R_{n-2})-R_{n-2} \, v'(R_{n-2})\bigr]\\
  =&~ R_{n-2} \, \bigl[ v_{\max}-v(R_{n-2})\bigr] - t \, \frac{R_{n-2}^2}{\ell_n} \, \bigl[ v_{\max}-v(R_{n-2})\bigr] \bigl[v_{\max}-v(R_{n-2})-R_{n-2} \, v'(R_{n-2})\bigr]\\
  \le&~ R_{n-2} \, \bigl[ v_{\max}-v(R_{n-2}) \bigr] \left[1-\frac{z_{n-2}}{\ell_n}\right].
\end{align*}
Since $z_{n-2}(0)=0$, a simple comparison argument shows that $z_{n-2}(t)\leq \ell_n$ for all times.

Next we prove that the inequality $z_{i+1}(t)\leq \ell_n$ being true for all $t\geq 0$ and for some $i\in\{1,\ldots,n-3\}$ implies $z_i(t) = t \, R_i(t) \, [ v(R_{i+1}(t)) - v(R_i(t)) ] \leq \ell_n$ for all $t\geq 0$.
We use the positive part $(z)_+ \doteq \max\{z,0\}$ and recall that $\sign_+(z_i) = \sign_+(v(R_{i+1})-v(R_i))= \sign_+(R_i-R_{i+1})$ for any $i\in\{1,\ldots,n-3\}$.
Let us compute
\begin{align*}
 \frac{d}{dt}(z_i)_+ =&~ R_i \, \bigl( v(R_{i+1})-v(R_i)\bigr)_+ + t \, \dot{R}_i \, \bigl( v(R_{i+1})-v(R_i)\bigr)_+ \\
  &~ + t \, R_i \, \bigl[ v'(R_{i+1}) \, \dot{R}_{i+1}- v'(R_{i}) \, \dot{R}_{i} \, \bigr] \, \sign_+\bigl(v(R_{i+1})-v(R_i)\bigr)\\
  =&~  R_i \bigl( v(R_{i+1})-v(R_i)\bigr)_+ \Biggl[1-\frac{(z_i)_+}{\ell_n}\Biggr] - v'(R_{i+1}) \, R_i \, R_{i+1} \frac{z_{i+1}}{\ell_n} \, \sign_+(z_i) + v'(R_i) \, R_i^2 \frac{(z_i)_+}{\ell_n}.
\end{align*}
The inequality $z_{i+1}\leq \ell_n$ and the assumption \eqref{V2} imply
\[%\begin{equation}\label{eq:oleinik_intermediate}
  \frac{d}{dt}(z_i)_+ \leq R_i \left[\vphantom{\frac{(z_i)_+}{\ell_n}} \, \bigl(v(R_{i+1})-v(R_i)\Bigr)_+ -v'(R_i) \, R_i\right] \left[1-\frac{(z_i)_+}{\ell_n}\right].
\]%\end{equation}
We observe that the first squared bracket on the right-hand-side of the above estimate is nonnegative.
Therefore, a comparison argument similar to that used before shows that $z_i(t)\leq \ell_n$ for all times $t\geq 0$.
Hence, the proof is complete.
\end{proof}

For $i\in\{1,\ldots,n-2\}$, the estimate \eqref{eq:oleinik} reads
\[\frac{v(R^n_{i+1}(t))-v(R^n_{i}(t))}{x^n_{i+1}(t)-x^n_i(t)}\leq \frac{1}{t},\]
which recalls the one-sided Lipschitz condition in \cite{oleinik} which characterises entropy solutions to \eqref{eq:CLintro}.

The result in Proposition~\ref{pro:oleinik} implies a uniform bound for $\rho^n$ in $\BVloc((0,+\infty)\times \R)$. In this sense, the $\L\infty\rightarrow\BV$ smoothing effect featured by genuinely nonlinear scalar conservation laws is intrinsically encoded in the particle scheme \eqref{eq:FTL}. In what follows, we denote by $\tv(f;\,U)$ the (local) total variation of a function $f$ on the subset $U\subset \R$.

\begin{proposition}\label{pro:bvloc}
%Assume $\bar\rho$ satisfies \eqref{I1} and $v$ satisfies \eqref{V1}, \eqref{V2}.
Assume that \eqref{I1}, \eqref{V1} and \eqref{V2} are satisfied.
Let $\delta>0$ and $a<b$. Then, the quantity
\[\sup_{t\ge\delta}\tv\bigl(\rho^n(t,\cdot);\,[a,b]\bigr)\]
is uniformly bounded with respect to $n$.
\end{proposition}

\begin{proof}
Fix $t \ge \delta$.
We assume that $x_0^n(t) \le a < b \le x_n^n(t)$, leaving to the reader the study of the remaining cases.
We introduce then
\begin{align*}
&I^n_a(t) \doteq \max\bigl\{i\in\{0,\ldots,n\}\colon x^n_i(t)\leq a \bigr\},
&I^n_b(t) \doteq \max\bigl\{i\in\{0,\ldots,n\}\colon x^n_i(t)\leq b \bigr\}.
\end{align*}
We consider
\begin{align*}
&\sigma^n(t,x) \doteq  v(\rho^n(t,x)) - \dfrac{1}{t} \, X^n(t),
&X^n(t,x) \doteq \sum_{i=0}^{n-1}x_i^n(t) ~ \mathbf{1}_{[x_i^n(t),x_{i+1}^n(t))}(x).
\end{align*}
We point out that $\sigma^n(t,\cdot)$ is non-increasing in $(x^n_0(t), x^n_n(t))$.
Indeed, by \eqref{eq:phantom}
\begin{align*}
&\sigma^n(t,x^n_i(t)^-) - \sigma^n(t,x^n_i(t)^+) = \dfrac{1}{t} [x^n_i(t) - x^n_{i-1}(t)] \ge 0,
&i \in \{1,n-1\},
\end{align*}
and the ODEs in \eqref{eq:FTL} together with the inequality \eqref{eq:oleinik} show that $\sigma^n(t,\cdot)$ is non-increasing in $(x^n_1(t), x^n_{n-1}(t))$.

By \eqref{eq:disdens} we can estimate the total variation of $v(\rho^n(t,\cdot))$ on $[a,b]$ as follows
\begin{align*}
  &~ \tv\Bigl(v(\rho^n(t,\cdot));\, [a,b]\Bigr)
  =\left|v(R^n_{I^n_a(t)+1}) - v(R^n_{I^n_a(t)})\right|
  +\tv\Bigl(v(\rho^n(t,\cdot));\, [x^n_{I^n_a(t)+1},x^n_{I^n_b(t)}]\Bigr)
  \\\leq&~
  \Bigl[v_{\max} - v(R)\Bigr]
  +\tv\Bigl(\sigma^n(t,\cdot);\, [x^n_{I^n_a(t)+1},x^n_{I^n_b(t)}]\Bigr)
  +\frac{1}{t} \, \tv\Bigl(X^n(t,\cdot);\, [x^n_{I^n_a(t)+1},x^n_{I^n_b(t)}]\Bigr)
  \\=&~
  \Bigl[v_{\max} - v(R)\Bigr]
  +\Bigl[\sigma^n(t,x^n_{I^n_a(t)+1})-\sigma^n(t,x^n_{I^n_b(t)})\Bigr]
  + \frac{1}{t} \, \bigl[X^n(t,x^n_{I^n_b(t)})-X^n(t,x^n_{I^n_a(t)+1})\Bigr]
  \\=&~
  \Bigl[v_{\max} - v(R)\Bigr]
  +\Bigl[v(\rho^n(t,x^n_{I^n_a(t)+1}))-v(\rho^n(t,x^n_{I^n_b(t)}))\Bigr]
  +\frac{2}{t} \, \Bigl[x^n_{I^n_b(t)} - x^n_{I^n_a(t)+1}\Bigr]
  \\\leq&~
  2\left[v_{\max} - v(R) +\frac{b-a}{\delta}\right].
\end{align*}
Since $v$ is monotone and continuous on $\R_+$, we get the assertion.
\end{proof}

Proposition~\ref{pro:bv_contraction} and Proposition~\ref{pro:bvloc} provide the needed compactness of $\rho^n$ with respect to the space variable.
Typically, in the context of scalar conservation laws (e.g.~the wave-front tracking scheme) %(see e.g.~the wave-front tracking scheme in \cite{bressan_book})
an $\L1$ uniform continuity estimate provides sufficient control of the time oscillations.
In our case, we are only able to provide a uniform time continuity estimate with respect to the \emph{$1$-Wasserstein distance}, which nevertheless will suffice to achieve strong $\L1$ compactness (with respect to both space and time).

We first recall the following concepts on the one dimensional $1$-Wasserstein distance.
Let $\mu$ be a probability measure on $\R$. We define the pseudo-inverse variable $X_\mu \in \L1([0,1])$ as
\[
X_\mu(z) \doteq \inf\{x\in\R\colon\mu((-\infty,x])>z\}.
\]
Given two probability measures $\mu$ and $\nu$ on $\R$, we set
\[W_1(\mu,\nu) \doteq \|X_{\mu}-X_{\nu}\|_{\L1([0,1])}.\]
By \eqref{eq:disdens} we have that
\[
X_{\rho^n(t,\cdot)}(z)
= \sum_{i=0}^{n-1} \, \Bigl[x^n_i(t) + \left(z-i\,\ell\right) R^n_i(t)^{-1}\Bigr] ~ \mathbf{1}_{[i\ell,(i+1) \, \ell)}(z).
\]

\begin{proposition}\label{pro:time_continuity}
Assume (I1) and (V1) are satisfied. There exists a constant $C$ independent of $n$, such that $W_1(\rho^n(t,\cdot),\rho^n(s,\cdot)) \leq C|t-s|$ for any $t,s>0$.
\end{proposition}

\begin{proof}
For $0<s<t$ we compute
\begin{align*}
  &~ W_1(\rho^n(t,\cdot),\rho^n(s,\cdot)) =  \|X_{\rho^n(t,\cdot)}-X_{\rho^n(s,\cdot)}\|_{\L1([0,1])}\\
  =&~ \sum_{i=0}^{n-1}\int_{i\ell}^{(i+1) \, \ell} \, \left|x^n_i(t)-x^n_i(s)+(z-i\,\ell)\left(R^n_i(t)^{-1} - R^n_i(s)^{-1}\right)\right|\d z\\
  \leq&~ \sum_{i=0}^{n-1} \ell \, |x^n_i(t)-x^n_i(s)| + \sum_{i=0}^{n-1}\left|R^n_i(t)^{-1} - R^n_i(s)^{-1}\right| \int_{i\ell}^{(i+1) \, \ell}(z-i\,\ell)\d z\\
  \leq&~ \max\{v_{\max},|v(R)|\} \, |t-s| + \sum_{i=0}^{n-1}\frac{\ell^2}{2}\int_{s}^{t}\left|\frac{d}{d\tau}\left(R^n_i(\tau)^{-1}\right)\right|\d \tau,
\end{align*}
and by using \eqref{eq:ODE_density} and \eqref{eq:phantom}
\begin{align*}
  &~  W_1(\rho^n(t,\cdot),\rho^n(s,\cdot))
  \\
  \leq&~
  \max\{v_{\max},|v(R)|\} \, |t-s|
  +\sum_{i=1}^{n-3} \ell \int_{s}^{t} |v(R^n_{i+1}(\tau))-v(R^n_i(\tau))| \d\tau
  +\ell \int_{s}^{t} |v_{\max} - v(R^n_{n-2}(\tau))| \d\tau
  \\
  \leq&~ \Bigl[\max\{v_{\max},|v(R)|\} + 2 [v_{\max} - v(R)] \Bigr] \, |t-s|.\qedhere
\end{align*}
\end{proof}

\begin{theorem}[Generalised Aubin-Lions lemma]\label{thm:aubin}
Let $T>0$, $a,b\in\R$ be fixed with $a<b$ and $v$ satisfy \eqref{V1}.
Let $\rho^n$ be a sequence in $\L\infty((0,T);\,\L1(\R))$ with $\rho^n(t,\cdot) \ge 0$ and $\|\rho^n(t,\cdot)\|_{\L1(\R)}=1$ for all $n \in \N$ and for all $t\in [0,T]$.
Assume further that
\begin{enumerate}[label={(\Alph*)},itemindent=*,leftmargin=0pt]\setlength{\itemsep}{0cm}
\item\label{A}
$\sup_{n\in \N} \left[ \int_0^T \left[\|v(\rho^n(t,\cdot))\|_{\L1([a,b])} + \tv(v(\rho^n(t,\cdot));\,[a,b])\right]\d t\right] < +\infty$,
\item\label{B}
there exists a constant $C>0$ independent of $n$ such that $W_1(\rho^n(t,\cdot),\rho^n(s,\cdot))\leq C|t-s|$ for all $s,t\in (0,T)$.
\end{enumerate}
Then, $\rho^n$ is strongly relatively compact in $\L1([0,T]\times [a,b])$.
\end{theorem}

The proof of Theorem~\ref{thm:aubin} is presented in the appendix~\ref{app}.

\begin{proof}[Conclusion of the proof of Theorem~\ref{thm:main}]

%Proposition~\ref{pro:bv_contraction} and Proposition~\ref{pro:bvloc} show that $\rho^n$ satisfies the assumption \ref{A} of Theorem~\ref{thm:aubin} on the time interval $[\delta,T]$ for arbitrary $0<\delta<T$ when $\bar\rho$ satisfies \eqref{I2} and $v$ satisfies \eqref{V2} respectively.
Proposition~\ref{pro:bv_contraction} and Proposition~\ref{pro:bvloc} show that $\rho^n$ satisfies the assumption \ref{A} of Theorem~\ref{thm:aubin} on the time interval $[\delta,T]$ for arbitrary $0<\delta<T$ when beside \eqref{I1} and \eqref{V1}, we assume either \eqref{I2} or \eqref{V2}.
The result in Proposition~\ref{pro:time_continuity} implies that $\rho^n$ satisfies assumption \ref{B} of Theorem~\ref{thm:aubin}.
Hence, by a simple diagonal argument stretching the time interval $[\delta,T]$ to $(0,T]$, one easily gets that $\rho^n$ has a subsequence (still denoted $\rho^n$) converging almost everywhere in $\Lloc1((0,T)\times \R)$.
Let $\rho$ be the limit of said subsequence.
\begin{enumerate}[label={\textbf{Step~\arabic*}},itemindent=*,leftmargin=0pt]\setlength{\itemsep}{0cm}%

\item\label{step1}\textbf{: $\boldsymbol\rho$ is a weak solution to \eqref{eq:cauchy}.}
Let $\varphi\in \Cc\infty(\R_+\times \R)$.
By \eqref{eq:disdens} we compute
\begin{align*}
  &~ \iint_{\R_+\times\R} \, \Bigl[ \rho^n(t,x) \, \varphi_t(t,x) + \rho^n(t,x) \, v(\rho^n(t,x)) \, \varphi_x(t,x) \Bigr] \d x \d t\\
  =&~ \sum_{i=0}^{n-1}\int_{\R_+} R^n_i(t)\left[\int_{x^n_i(t)}^{x^n_{i+1}(t)} \varphi_t(t,x)  \d x
  +v(R^n_i(t)) \Bigl[ \varphi(t,x_{i+1}^n(t)) - \varphi(t,x_{i}^n(t))\Bigr] \right]\d t\\
  =&~ \sum_{i=0}^{n-1}\int_{\R_+} R^n_i(t)
  \Biggl[ \frac{d}{dt}\left(\int_{x^n_i(t)}^{x^n_{i+1}(t)}\varphi(t,x)\d x\right)
  +\Bigl[ \dot{x}_i^n(t) - v(R^n_i(t)) \Bigr] \varphi(t,x^n_i(t))
  \\
  &~ \hphantom{\sum_{i=0}^{n-1}\int_{\R_+} R^n_i(t) \Biggl[}
  -\Bigl[ \dot{x}^n_{i+1}(t) - v(R^n_i(t)) \Bigr] \varphi(t,x^n_{i+1}(t))
  \Biggr] \d t
  \\
  =&~ \sum_{i=0}^{n-1}\int_{\R_+} \, \Biggl[ -\dot{R}^n_i(t)\left(\int_{x^n_i(t)}^{x^n_{i+1}(t)} \varphi(t,x) \d x\right)
  + R^n_i(t) \Bigl[ \dot{x}_i^n(t)-v(R^n_i(t)) \Bigr] \varphi(t,x^n_i(t))
  \\
  &~ \hphantom{\sum_{i=0}^{n-1}\int_{\R_+} \, \Biggl[}
  - \dfrac{R^n_i(t)^2}{\ell} \, \Bigl[ \dot{x}^n_{i+1}(t)-v(R^n_i(t)) \Bigr] \Biggl[\int_{x^n_i(t)}^{x^n_{i+1}(t)} \varphi(t,x^n_{i+1}(t)) \d x\Biggr]
  \Biggr]\d t
  -\int_{\R} \rho^n(0,x) \, \varphi(0,x) \d x.
\end{align*}
By \eqref{eq:mass} and the definition of $R^n_i$ we have that
\begin{align*}
&~\modulo{\int_{\R} \, \Bigl[\bar{\rho}(x) - \rho^n(0,x)\Bigr] \varphi(0,x) \d x}
\\\le&~
\int_{-\infty}^{\bar{x}^n_{0}} \bar{\rho}(x) \, \varphi(0,x) \d x
+\int_{\bar{x}^n_n}^{+\infty} \bar{\rho}(x) \, \varphi(0,x) \d x
+\sum_{i=0}^{n-1} \modulo{\int_{\bar{x}^n_i}^{\bar{x}^n_{i+1}} \, \Bigl[\bar{\rho}(x) - R^n_i(0)\Bigr] \varphi(0,x) \d x}
\\\le&~
2\ell_n \, \|\varphi(0,\cdot)\|_{\L\infty(\R)}
+\sum_{i=0}^{n-1} \modulo{\int_{\bar{x}^n_i}^{\bar{x}^n_{i+1}} \bar{\rho}(x) \left[ \varphi(0,x) - \fint_{\bar{x}^n_i}^{\bar{x}^n_{i+1}} \varphi(0,y) \d y
\right] \d x}
\end{align*}
and clearly the above quantity goes to zero as $n \to + \infty$.
Now we have to consider two separate cases.

\begin{enumerate}[label={\sc Case~\arabic*},itemindent=*,leftmargin=0pt]\setlength{\itemsep}{0cm}%

\item\label{case1}\textsc{:~$\bar\rho$ is compactly supported.}
In this case, we can use the improved construction of the particle scheme described in Remark~\ref{rem01} and the equations analogous to \eqref{eq:ODE_density} and \eqref{eq:FTL} as follows.
Assuming that $\supp[\varphi]\subset [\delta,T]\times \R$ for some $0<\delta<T$, we obtain
\begin{align*}
&~
\Biggl|\iint_{\R_+\times\R} \, \bigl[\rho^n(t,x) \, \varphi_t(t,x) + \rho^n(t,x) \, v(\rho^n(t,x)) \, \varphi_x(t,x)\bigr] \d x \d t \Biggr|
\\
=&~ \Biggl| \sum_{i=0}^{n-2}\int_0^T \frac{R^n_i(t)^2}{\ell} \, \bigl[v(R^n_{i+1}(t)) - v(R^n_i(t))\bigr]
\left[\int_{x^n_i(t)}^{x^n_{i+1}(t)}
\bigl[\varphi(t,x)-\varphi(t,x^n_{i+1}(t))\bigr]  \d x\right] \d t
\\
&~ + \int_0^T \frac{R^n_{n-1}(t)^2}{\ell} \, \bigl[v_{\max}-v(R^n_{n-1}(t))\bigr]
\left[\int_{x^n_{n-1}(t)}^{x^n_{n}(t)}
\bigl[\varphi(t,x)-\varphi(t,x^n_{n}(t))\bigr]  \d x\right] \d t \Biggr|
\\
\leq&~
\frac{T \, \lip[\varphi] \, \ell}{2}\sup_{t\in [\delta,T]} \left[\sum_{i=0}^{n-2} \, \bigl|v(R^n_{i+1}(t))-v(R^n_i(t))\bigr| + \bigl|v_{\max}-v(R^n_{n-1}(t))\bigr| \right]
\\
\leq&~
\frac{T \, \lip[\varphi] \, \ell}{2}
\left[v_{\max} - v(R)+\sup_{t\in [\delta,T]} \tv\bigl(v(\rho^n(t,\cdot));\,J(T)\bigr)\right],
\label{eq:weak_compact}\tag{$\spadesuit$}
\end{align*}
where $J(T) \doteq \bigl[\min\{\supp[\bar\rho]\} + v(R) \, T , \max\{\supp[\bar\rho]\} + v_{\max} \, T \bigr]$.
%where $J(T) \doteq \bigl[\min\{\overline{\supp[\bar\rho]}\} + v(R) \, T , \max\{\overline{\supp[\bar\rho]}\} + v_{\max} \, T \bigr]$.
Hence, by Proposition~\ref{pro:bvloc} the right hand side in \eqref{eq:weak_compact} tends to zero as $n\rightarrow +\infty$, and since $\rho^n$ tends to $\rho$ almost everywhere up to a subsequence we have that $\rho$ is a weak solution to the Cauchy problem \eqref{eq:cauchy} for positive times.

\item\label{case2}\textsc{:~$\bar\rho$ is NOT compactly supported.}
For simplicity we shall assume that $\supp[\bar\rho]$ is unbounded both from above and from below.
The remaining cases are minor variations of this one.
Assume $\supp[\varphi]\subset [\delta,T]\times [a,b]$ for some $0<\delta<T$ and for some $a<b$.
Let $n\in\N$ be sufficiently large so that $\bar{x}^n_1 < a-v_{\max} \, T$ and $\bar{x}^n_{n-1} > b-v(R) \, T$.
Such a choice is possible because $\supp[\bar\rho]$ is unbounded both from above and from below, which implies that the sequence $\supp[\rho^n(0,\cdot)]$ is not uniformly bounded with respect to $n\in\N$ both from above and from below.
Such assumptions imply that $x_1^n(t)<a$ and $x^n_{n-1}(t)>b$ for all $t\in[0,T]$. We have
\begin{align*}
  &~ \Biggl|\iint_{\R_+\times\R} \left[\rho^n(t,x) \, \varphi_t(t,x) + \rho^n(t,x) \,  v(\rho^n(t,x)) \, \varphi_x(t,x)\right]\d x \d t \Biggr|\\
  =&~ \Biggl|\sum_{i=1}^{n-2} \int_{\R_+} R^n_i(t)\left[\int_{x^n_i(t)}^{x^n_{i+1}(t)} \varphi_t(t,x)  \d x
  +v(R^n_i(t)) \Bigl[ \varphi(t,x_{i+1}^n(t)) - \varphi(t,x_{i}^n(t))\Bigr] \right]\d t \Biggr|
\end{align*}
for all $\varphi\in \Cc\infty(\R_+\times \R)$ and the assertion can be obtained as in \ref{case1} (we omit the details).
\end{enumerate}

\item\textbf{: $\boldsymbol\rho$ satisfies the entropy inequality in Definition~\ref{def:entropy_sol}.}
Let $\varphi\in \Cc\infty((0,+\infty) \times \R)$ with $\varphi\geq 0$ and $k\geq 0$.
By \eqref{eq:disdens}
\begin{align*}
&~\iint_{\R_+\times\R} \, \biggl[|\rho(t,x)-k| \, \varphi_t(t,x) + \sign(\rho(t,x) - k) \Bigl[ f(\rho(t,x)) - f(k)\Bigr] \varphi_x(t,x)\biggr] \d x\d t
\\
=&~
\int_{\R_+} \int_{-\infty}^{x_0^n(t)} \, \biggl[k \, \varphi_t(t,x) + f(k) \, \varphi_x(t,x)\biggr] \d x\d t
+\int_{\R_+} \int^{+\infty}_{x_n^n(t)} \, \biggl[k \, \varphi_t(t,x) + f(k) \, \varphi_x(t,x)\biggr] \d x\d t
\\&~
+\sum_{i=0}^{n-1} \int_{\R_+}
\Biggl[|R^n_i(t)-k| \left( \int_{x^n_{i}(t)}^{x_{i+1}^n(t)} \varphi_t(t,x) \d x \right)
\\
&~\hphantom{+\sum_{i=0}^{n-1} \int_{\R_+}\Biggl[}
+ \sign(R^n_i(t) - k) \Bigl[ f(R^n_i(t)) - f(k)\Bigr] \Bigl[\varphi(t,x_{i+1}^n(t)) - \varphi(t,x_{i}^n(t)) \Bigr] \Biggr] \d t
\\
=&~
k \int_{\R_+} \, \Bigl[
\bigl[ v(k) - \dot{x}_0^n(t) \bigr] \varphi(t,x_0^n(t))
-
\bigl[ v(k) - \dot{x}_n^n(t) \bigr] \varphi(t,x_n^n(t))
\Bigr] \d t
\\&~
+\sum_{i=0}^{n-1} \int_{\R_+}
\sign(R^n_i(t) - k)
\Biggl[ \Bigl[R^n_i(t)-k\Bigr] \, \dfrac{d}{dt} \left( \int_{x^n_{i}(t)}^{x_{i+1}^n(t)} \varphi(t,x) \d x \right)
\\
&~\hphantom{+\sum_{i=0}^{n-1} \int_{\R_+}\Bigl[ \sign(R^n_i(t) - k)\Biggl[}
+ \Bigl[ f(R^n_i(t)) - f(k) - (R^n_i(t) - k) \, \dot{x}_{i+1}^n(t) \Bigr] \varphi(t,x_{i+1}^n(t))
\\
&~\hphantom{+\sum_{i=0}^{n-1} \int_{\R_+}\Bigl[ \sign(R^n_i(t) - k)\Biggl[}
- \Bigl[ f(R^n_i(t)) - f(k) - (R^n_i(t) - k) \, \dot{x}_{i}^n(t) \Bigr]  \varphi(t,x_{i}^n(t)) \Biggr] \d t
\\
=&~
k \int_{\R_+} \, \Bigl[
\bigl[ v(k) - \dot{x}_0^n(t) \bigr] \varphi(t,x_0^n(t))
-
\bigl[ v(k) - \dot{x}_n^n(t) \bigr] \varphi(t,x_n^n(t))
\Bigr] \d t
\\&~
+\sum_{i=0}^{n-1} \int_{\R_+}
\sign(R^n_i(t) - k)
\Biggl[ -\dot{R}^n_i(t) \, \left( \int_{x^n_{i}(t)}^{x_{i+1}^n(t)} \varphi(t,x) \d x \right)
\\
&~\hphantom{+\sum_{i=0}^{n-1} \int_{\R_+}\Bigl[ \sign(R^n_i(t) - k)\Biggl[}
- \Bigl[ R^n_i(t) \bigl[ \dot{x}_{i+1}^n(t) - v(R^n_i(t)) \bigr] - k \bigl[ \dot{x}_{i+1}^n(t) -  v(k) \bigr] \Bigr] \varphi(t,x_{i+1}^n(t))
\\
&~\hphantom{+\sum_{i=0}^{n-1} \int_{\R_+}\Bigl[ \sign(R^n_i(t) - k)\Biggl[}
+ \Bigl[R^n_i(t) \bigl[ \dot{x}_{i}^n(t) - v(R^n_i(t)) \bigr] - k \bigl[ \dot{x}_{i}^n(t) - v(k) \bigr] \Bigr]  \varphi(t,x_{i}^n(t)) \Biggr] \d t.
\end{align*}
Now we have to consider two separate cases.

\begin{enumerate}[label={\sc Case~\arabic*},itemindent=*,leftmargin=0pt]\setlength{\itemsep}{0cm}%

\item\label{case21}\textsc{:~$\bar\rho$ is compactly supported.}
In this case, we can use the improved construction of the particle scheme described in Remark~\ref{rem01} and the equations analogous to \eqref{eq:ODE_density} and \eqref{eq:FTL} as follows.
Assuming that $\supp[\varphi]\subset [\delta,T]\times \R$ for some $0<\delta<T$, we obtain
\begin{align*}
&\iint_{\R_+\times\R} \, \Biggl[|\rho(t,x)-k| \, \varphi_t(t,x) + \sign(\rho(t,x) - k) \Bigl[ f(\rho(t,x)) - f(k)\Bigr] \varphi_x(t,x)\Biggr] \d x\d t
\\
=&
k \int_{\R_+} \, \Bigl[
\bigl[ v(k) - v(R^n_0(t)) \bigr] \varphi(t,x_0^n(t))
-
\bigl[ v(k) - v_{\max} \, \bigr] \varphi(t,x_n^n(t))
\Bigr] \d t
\\&
+\sum_{i=0}^{n-2} \int_{\R_+}\!
\sign(R^n_i(t) - k)
\Biggl[ \dfrac{R^n_i(t)^2}{\ell_n} \Bigl[v(R^n_{i+1}(t))-v(R^n_i(t))\Bigr] \Biggl[ \int_{x^n_{i}(t)}^{x_{i+1}^n(t)} \bigl[ \varphi(t,x) - \varphi(t,x_{i+1}^n(t)) \bigr]\d x \Biggr]
\\
&\hphantom{+\sum_{i=0}^{n-2} \int_{\R_+}}
+ k \Bigl[ \bigl[ v(R_{i+1}^n(t)) - v(k) \bigr] \varphi(t,x_{i+1}^n(t))
- \bigl[ v(R_{i}^n(t)) - v(k) \bigr]  \varphi(t,x_{i}^n(t)) \Bigr] \Biggr] \d t
\\
&
+\int_{\R_+}
\sign(R^n_{n-1}(t) - k)
\Biggl[ \dfrac{R^n_{n-1}(t)^2}{\ell_n} \, \Bigl[v_{\max}-v(R^n_{n-1}(t))\Bigr] \Biggl[ \int_{x^n_{n-1}(t)}^{x_{n}^n(t)} \, \bigl[ \varphi(t,x) - \varphi(t,x_{n}^n(t)) \bigr] \d x \Biggr]
\\
&\hphantom{+\int_{\R_+}}
+k \Bigl[ \bigl[v_{\max} - v(k) \bigr] \varphi(t,x_{n}^n(t))
- \bigl[ v(R_{n-1}^n(t)) - v(k) \bigr] \varphi(t,x_{n-1}^n(t)) \Bigr] \Biggr] \d t.
\end{align*}
We already proved, see \eqref{eq:weak_compact}, that
\begin{align*}
&
\sum_{i=0}^{n-2} \int_{\R_+}
\sign(R^n_i(t) - k)
\dfrac{R^n_i(t)^2}{\ell_n} \, \Bigl[v(R^n_{i+1}(t))-v(R^n_i(t))\Bigr] \Biggl[ \int_{x^n_{i}(t)}^{x_{i+1}^n(t)} \, \bigl[ \varphi(t,x) - \varphi(t,x_{i+1}^n(t)) \bigr]\d x \Biggr] \d t
\\
&\hphantom{\Biggl|}
+\int_{\R_+}
\sign(R^n_{n-1}(t) - k)
\dfrac{R^n_{n-1}(t)^2}{\ell_n} \, \Bigl[v_{\max}-v(R^n_{n-1}(t))\Bigr] \Biggl[ \int_{x^n_{n-1}(t)}^{x_{n}^n(t)} \, \bigl[ \varphi(t,x) - \varphi(t,x_{n}^n(t)) \bigr] \d x \Biggr] \d t
\end{align*}
converges to zero as $n\to+\infty$.
Hence, to conclude it suffices to observe that
\begin{align*}
&~ k \Biggl[
\bigl[ v(k) - v(R^n_0(t)) \bigr] \varphi(t,x_0^n(t))
-
\bigl[ v(k) - v_{\max} \, \bigr] \varphi(t,x_n^n(t))
\\
&~
+\sum_{i=0}^{n-2}
\sign(R^n_i(t) - k) \,
\Bigl[ \bigl[ v(R_{i+1}^n(t)) - v(k) \bigr] \varphi(t,x_{i+1}^n(t))
- \bigl[ v(R_{i}^n(t)) - v(k) \bigr]  \varphi(t,x_{i}^n(t)) \Bigr]
\\
&~\hphantom{+k \Biggl[}
+\sign(R^n_{n-1}(t) - k)
\Bigl[ \bigl[v_{\max} - v(k) \bigr] \varphi(t,x_{n}^n(t))
- \bigl[ v(R_{n-1}^n(t)) - v(k) \bigr] \varphi(t,x_{n-1}^n(t)) \Bigr] \Biggr]
\\
=&~
k \Biggl[ \sum_{i=1}^{n-1}
\bigl[ \sign(R^n_{i-1}(t) - k) - \sign(R^n_i(t) - k) \bigr] \bigl[ v(R_{i}^n(t)) - v(k) \bigr] \varphi(t,x_{i}^n(t))
\\
&~\hphantom{k \Biggl[}
+ \bigl[ 1 +\sign(R^n_0(t) - k) \bigr] \bigl[ v(k) - v(R^n_0(t)) \bigr] \varphi(t,x_{0}^n(t))
\\
&~\hphantom{k \Biggl[}
+\bigl[1+\sign(R^n_{n-1}(t) - k)\bigr] \bigl[v_{\max} - v(k) \bigr] \varphi(t,x_{n}^n(t)) \Biggr]
\ge 0.
\end{align*}
\item\textsc{:~$\bar\rho$ is NOT compactly supported.}
For simplicity we shall assume that $\supp[\bar\rho]$ is unbounded both from above and from below.
The remaining cases are minor variations of this one.
%Assume $\supp[\varphi]\subset [\delta,T]\times [a,b]$ for some $0<\delta<T$ and for some $a,b\in \R$ with $a<b$.
%Let $n\in\N$ be sufficiently large so that $\bar{x}^n_1 < a-v_{\max} \, T$ and $\bar{x}^n_{n-1} > b-v(R) \, T$.
%Such a choice is always possible because $\supp[\bar\rho]$ is unbounded both from above and from below, which implies that the sequence $\supp[\rho^n(0,\cdot)]$ is not uniformly bounded with respect to $n\in\N$ both from above and from below.
%Such assumptions imply that $x_1^n(t)<a$ and $x^n_{n-1}(t)>b$ for all $t\in[0,T]$. We have
Then, with the same notations and assumptions used in \ref{case2} of \ref{step1}, we have

\begin{align*}
&~\iint_{\R_+\times\R} \, \biggl[|\rho(t,x)-k| \, \varphi_t(t,x) + \sign(\rho(t,x) - k) \Bigl[ f(\rho(t,x)) - f(k)\Bigr] \varphi_x(t,x)\Biggr] \d x\d t
\\
=&~
\sum_{i=1}^{n-2} \int_{\R_+}
\Biggl[|R^n_i(t)-k| \left( \int_{x^n_{i}(t)}^{x_{i+1}^n(t)} \varphi_t(t,x) \d x \right)
\\
&~\hphantom{\sum_{i=0}^{n-1} \int_{\R_+}\Biggl[}
+ \sign(R^n_i(t) - k) \Bigl[ f(R^n_i(t)) - f(k)\Bigr] \Bigl[\varphi(t,x_{i+1}^n(t)) - \varphi(t,x_{i}^n(t)) \Bigr] \Biggr] \d t
\end{align*}
for all $\varphi\in \Cc\infty((0,+\infty) \times \R)$ and the assertion can be obtained as in the above \ref{case21} (we omit the details).
\qedhere
\end{enumerate}
\end{enumerate}

\end{proof}

\section{Numerical simulations}

\begin{figure}[!ht]
%\begin{minipage}[l]{.4\textwidth}
%\begin{center}
%\includegraphics[width=\textwidth]{figures/LWRC_0408_400_253}
%\end{center}
%\end{minipage}
\begin{center}
\includegraphics[width=.33\textwidth]{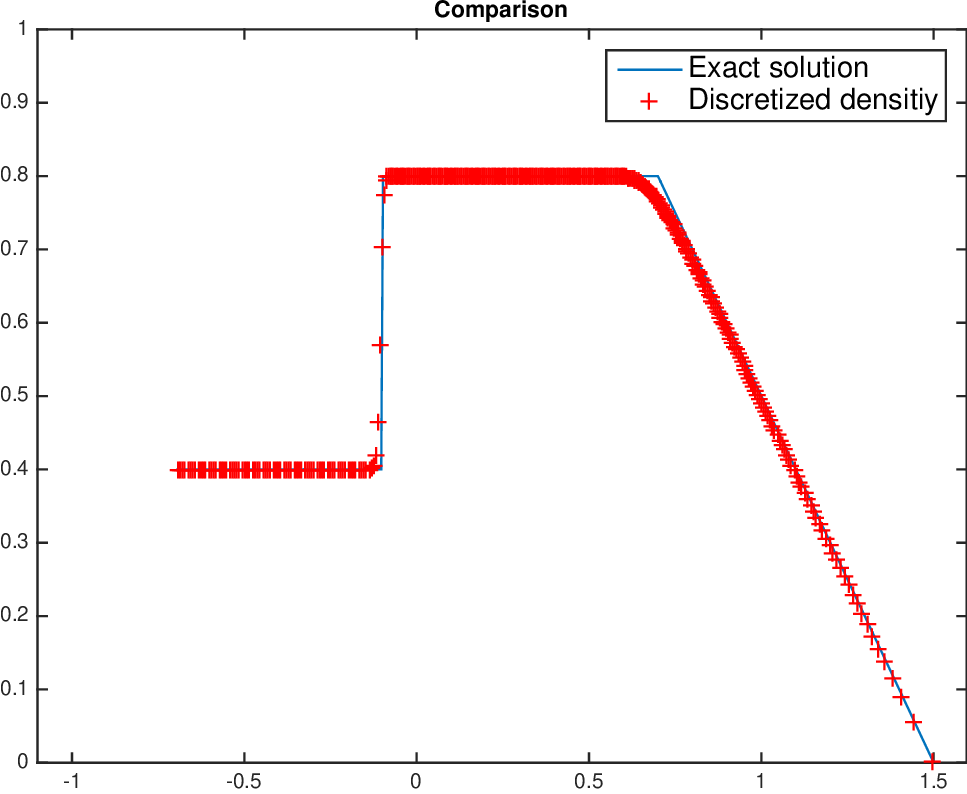}
\end{center}
\caption{Comparison between the exact solution (continuous blue line in the pdf version of the paper) and $\rho^n(t,x)$ (``$+$'' in red in the pdf version of the paper) for $N=400$ particles and initial datum \eqref{IC}. \label{fig:confr}}
\end{figure}

This section is devoted to present numerical simulations for the particle method described above.
We compare the numerical simulations with the exact solutions obtained by the method of characteristics.
\begin{figure}[!ht]
%\begin{minipage}[l]{.4\textwidth}
%\begin{center}
%\includegraphics[width=\textwidth]{figures/LWR_0408_200_1}
%\end{center}
%\end{minipage}\\
%\begin{minipage}[c]{.4\textwidth}
%\begin{center}
%\includegraphics[width=\textwidth]{figures/LWR_0408_200_50}
%\end{center}
%\end{minipage}\\
%\begin{minipage}[c]{.4\textwidth}
%\begin{center}
%\includegraphics[width=\textwidth]{figures/LWR_0408_200_101}
%\end{center}
%\end{minipage}
\begin{center}
\includegraphics[width=.32\textwidth]{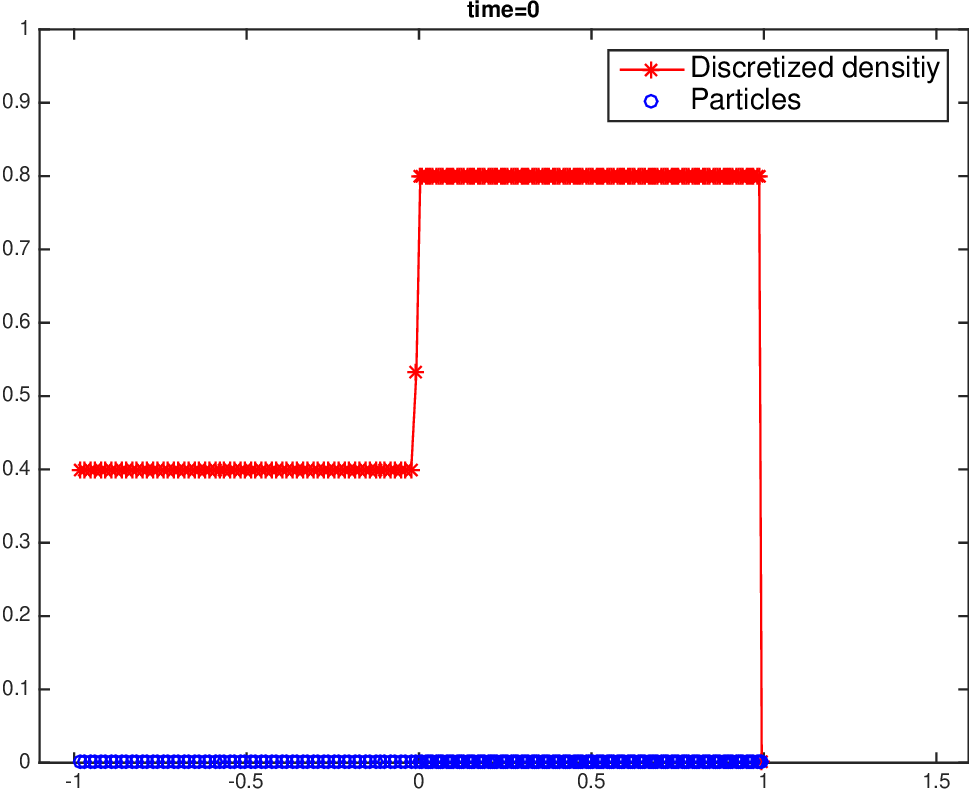}~
\includegraphics[width=.32\textwidth]{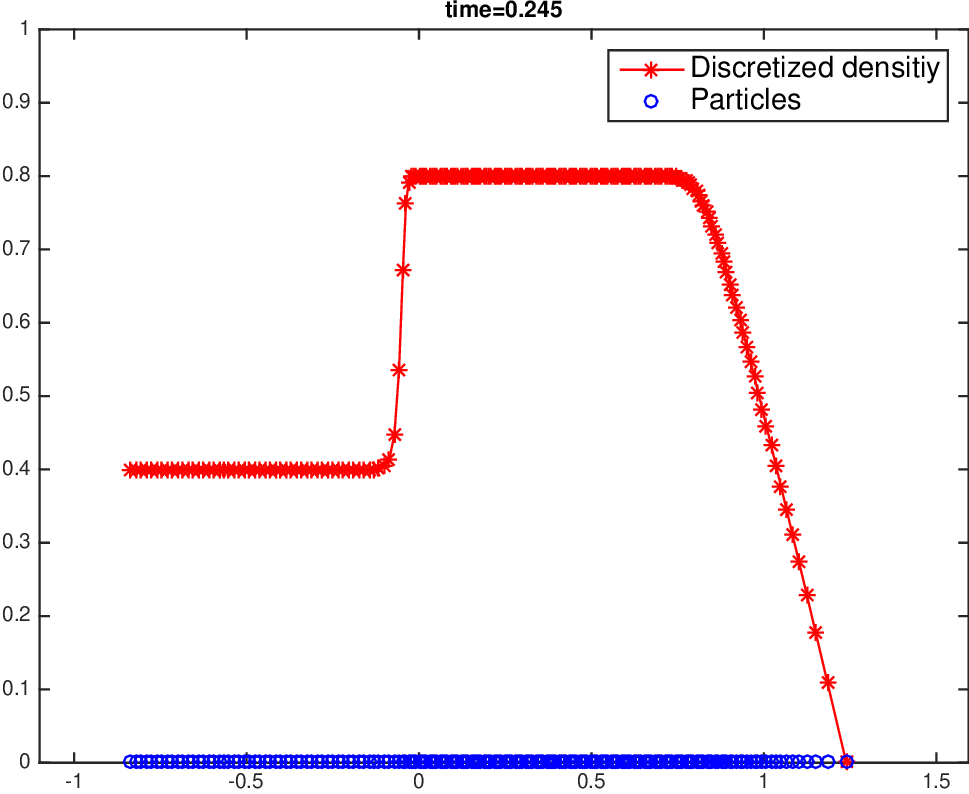}~
\includegraphics[width=.32\textwidth]{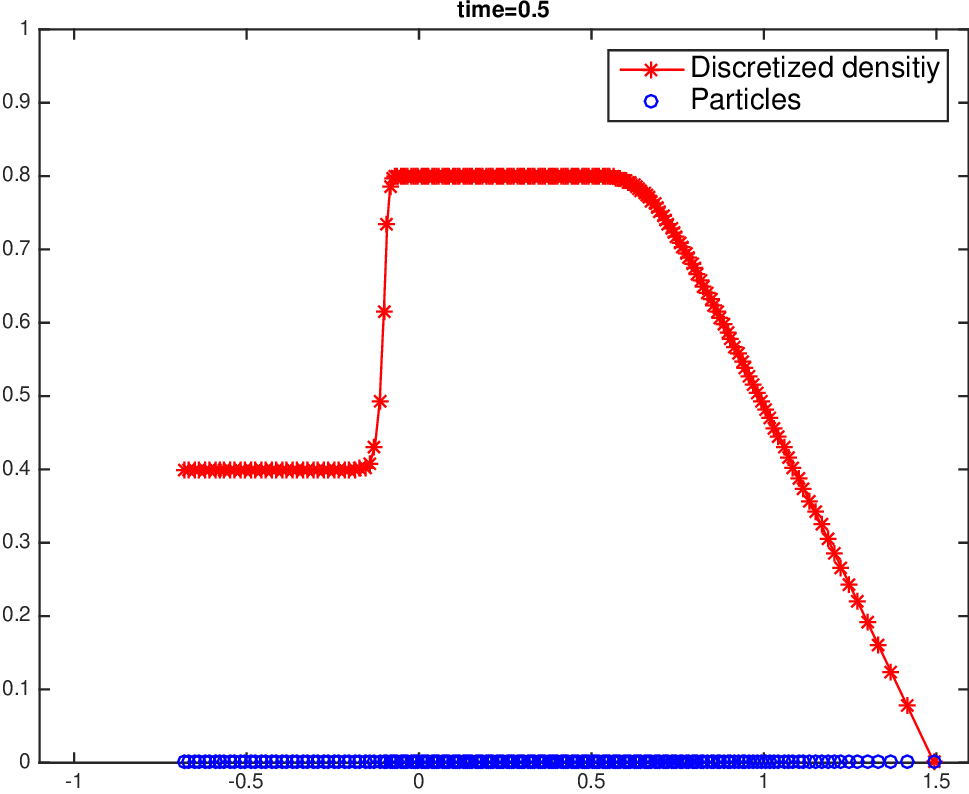}
\end{center}
\caption{The evolution of $\rho^n$ with  initial datum \eqref{IC}.
The cirles in the bottom (in blue in the pdf version of the paper) denote particle location, while the stars in the top (in red in the pdf version of the paper) denote the computed density.}
\label{fig:Test1}
\end{figure}
The particle system \eqref{eq:FTL} is solved using the Runge-Kutta MATLAB solver ODE23, with the  initial mesh size determined by the total number of particles $N$ and the initial density values.
In \figurename~\ref{fig:Test1} we take $N=200$ particles and the initial datum
\begin{equation}\label{IC}
\bar{\rho}(x)=\begin{cases}
 0.4 &\text{if } -1\leq x\leq 0, \\
 0.8 &\text{if } 0< x\leq 1, \\
 0 & \text{otherwise},
\end{cases}
\end{equation}
and final time $t=0.5$.
In \figurename~\ref{fig:confr} we compare the simulation with $N=400$ particles with exact solutions and final time $t=0.5$.% built using characteristics method.

For several values of $N$, we do a quantitative evaluation through the discrete $\L1$-error, computed as the difference between approximated and exact solutions.
The results are collected in Table~\ref{errors}.
\begin{table}[!ht]
\begin{center}
\begin{tabular}{|c|c|c|c|c|c|}
\hline
 $N$\vphantom{$\displaystyle\sum$} & 50 & 100 & 200 & 400 & 1000\\
 \hline
 \hline
  \vphantom{$\displaystyle\sum$}  & $4.8e-02$ & $2.9e-02$ & $1.4e-02$ & $8.2e-03$ & $3.6e-03$\\
 \hline
\end{tabular}
\end{center}
 \caption{Discrete $\L1$-errors corresponding to different numbers of particles $N$.\label{errors}}
\end{table}

%\section*{Acknowledgments}
%
%The authors would like to thanks Giovanni Russo for comments and suggestions on the numerical part.

\appendix

\section{Proof of Theorem~\ref{thm:aubin}}\label{app}

We recall the following theorem.

\begin{theorem}[\cite{rossisavare}]\label{thm:rs}
    Let $\mathbb{X}$ be a separable Banach space. Let
    \begin{enumerate}[label=(\Alph*),itemindent=*,leftmargin=0pt]\setlength{\itemsep}{0cm}%
    \setcounter{enumi}{5}
        \item\label{F} $\mathfrak{F} \,\colon\, \mathbb{X}\to[0,+\infty]$ be a \emph{normal coercive integrand}, i.e.\ $\mathfrak{F}$ is lower semi-continuous w.r.t.\ the topology of $\mathbb{X}$ and its sub-levels are relatively compact in $\mathbb{X}$;
    \end{enumerate}
    \begin{enumerate}[label=(\alph*),itemindent=*,leftmargin=0pt]\setlength{\itemsep}{0cm}%
    \setcounter{enumi}{6}
        \item\label{g} $\mathfrak{g} \,\colon\, \mathbb{X}\times \mathbb{X}\to[0,+\infty]$ be a \emph{pseudo-distance}, i.e.
        $\mathfrak{g}$ is lower semi-continuous w.r.t. the topology of $\mathbb{X}$,
        and if $\nu,\mu\in \mathbb{X}$ are such that $\mathfrak{g}(\nu,\mu)=0$, $\mathfrak{F}[\nu]<+\infty$ and $\mathfrak{F}[\mu]<+\infty$, then $\nu=\mu$.
    \end{enumerate}
    For a fixed $T>0$, let $U$ be a set of measurable functions $\nu \colon (0,T) \to \mathbb{X}$ such that
    \begin{align}
        \label{eq:savare_hypo}
        &\sup_{\nu\in U}\int_0^T \mathfrak{F}\left[\nu(t)\right] \d t<+\infty&
        &\text{and}&
        \lim_{h \searrow 0} \left[\sup_{\nu\in U}\int_0^{T-h} \mathfrak{g}\left(\nu(t+h),\nu(t)\right) \d t\right]=0.
    \end{align}
    Then $U$ is strongly relatively compact in $\L1((0,T); \mathbb{X})$.
\end{theorem}

Let $I \doteq [a,b]$.
With the same notation of Theorem~\ref{thm:rs}, we set $\mathbb{X} \doteq \L1(I)$, $U \doteq \{\rho^n\}_n$, and
\[\mathfrak{F}[\rho] \doteq \|v(\rho)\|_{\L1(I)} + \tv(v(\rho);\,I).\]
Given a probability measure $\mu$, we set
\[\tilde{\mu} \doteq \mu|_{(a,b)} + \mu((-\infty,a]) \, \delta_a + \mu([b,+\infty)) \, \delta_b.\]
We then define
\[\mathfrak{g}(\mu,\nu) \doteq \begin{cases}
W_1(\tilde{\mu},\tilde{\nu}) & \hbox{if } \mu(\R)=\nu(\R)=1, \\
+\infty & \hbox{otherwise}.
\end{cases}\]

The lower semi-continuity of $\mathfrak{F}$ with respect to $\L1(I)$ follows from \cite[Theorem 1, page 172]{EvansG} and from the fact that $v$ is continuous.
The compactness property follows from \cite[Theorem 4, page 176]{EvansG}.
This proves that $\mathfrak{F}$ satisfies the assumption \ref{F}.
Let $\mu,\nu \in \L1(I)$ be two probability measures.
We observe that $W_1(\tilde{\mu},\tilde{\nu})=\|X_{\tilde{\mu}}-X_{\tilde{\nu}}\|_{\L1([0,1])}$, with
\[X_{\tilde{\mu}} \doteq a \, \mathbf{1}_{[0,\mu((-\infty,a])]} + X_{\mu} \, \mathbf{1}_{(\mu((-\infty,a]),\mu([b,+\infty)))} + b \, \mathbf{1}_{[\mu([b,+\infty)),1]}.\]
Consequently, setting $F_{\tilde{\mu}},F_{\tilde{\nu}} \colon I\rightarrow [0,1]$
\begin{align*}
&F_{\tilde{\mu}}(x) \doteq \int_{-\infty}^x \mu(y)\d y,
&F_{\tilde{\nu}}(x)\doteq \int_{-\infty}^x \nu(y)\d y,
\end{align*}
we easily get, from the fundamental theorem of integral calculus,
\[
W_1(\tilde{\mu},\tilde{\nu})
=\|X_{\tilde{\mu}}-X_{\tilde{\nu}}\|_{\L1([0,1])}
=\int_I |F_{\tilde{\mu}}(x)-F_{\tilde{\nu}}(x)|\d x
\leq \int_a^b \int_a^x |\mu(y)-\nu(y)|\d y \d x \leq (b-a)\|\mu-\nu\|_{\L1([a,b])},
\]
and this implies the (lower semi) continuity of $\mathfrak{g}$ with respect to $\L1(I)$.
The remaining part of the assumption \ref{g} is straightforward.
Finally, the conditions \eqref{eq:savare_hypo} easily follow from \ref{A} and \ref{B} in the statement of Theorem~\ref{thm:aubin}.

\section*{Acknowledgments}

MDF is supported by the Italian MIUR-PRIN project $2012L5WXHJ\_003$.
SF is partially supported by the Italian INdAM-GNAMPA 2015 mini-project: Analisi e stabilit\`a per modelli di equazioni alle derivate parziali nella matematica applicata. MDF acknowledges the Gran Sasso Science Institute in L'Aquila for the opportunity to teach this topic in a mini-course in March 2016. Part of these notes is based on the outcome of said mini-course. The authors acknowledge useful suggestions by Giovanni Russo for the numerical part.

%\begin{acknowledgements}
%The authors would like to thanks Giovanni Russo for comments and suggestions on the numerical part.
%\end{acknowledgements}


%\bibliographystyle{plain}
%\bibliography{UMI}

\end{document}